\documentclass[11pt,fleqn,twoside]{article}

\usepackage{amsmath}
\usepackage{amsthm}
\usepackage{amssymb}
\makeatletter
\newcommand{\prava}[1]{\small\it
\begin{flushleft}
Copyright \copyright \ 2000 by  #1
\end{flushleft}}

\newcommand{\name}[1]{\begin{flushleft}
                       \LARGE \bf #1
                       \end{flushleft}\vspace{-3mm}}

\newcommand{\Author}[1]{\begin{flushleft}
                       \it #1 \end{flushleft}}

\newcommand{\Adress}[1]{\begin{flushleft}
                       \it #1 \end{flushleft}}

\newcommand{\Date}[1]{\begin{flushleft}
                      \small  \it #1 \end{flushleft}}

\newcommand{\ehkol}{Author \ name}
\newcommand{\ohkol}{Article \ name}

\renewcommand{\@evenhead}{
\hspace*{-3pt}\raisebox{-15pt}[\headheight][0pt]{\vbox{\hbox to \textwidth 
{\thepage \hfil \ehkol}\vskip4pt \hrule}}}
\renewcommand{\@oddhead}{
\hspace*{-3pt}\raisebox{-15pt}[\headheight][0pt]{\vbox{\hbox to \textwidth 
{\ohkol \hfil \thepage}\vskip4pt\hrule}}}
\renewcommand{\@evenfoot}{}
\renewcommand{\@oddfoot}{}

\long\def\@makecaption#1#2{%
  \vskip\abovecaptionskip
  \sbox\@tempboxa{\small \textbf{#1.}\ \ #2}%
  \ifdim \wd\@tempboxa >\hsize
    {\small \textbf{#1.}\ \ #2}\par
  \else
    \global \@minipagefalse
    \hb@xt@\hsize{\hfil\box\@tempboxa\hfil}%
  \fi
  \vskip\belowcaptionskip}

\def\numberwithin#1#2{\@ifundefined{c@#1}{\@nocounterr{#1}}{%
  \@ifundefined{c@#2}{\@nocnterr{#2}}{%
  \@addtoreset{#1}{#2}%
  \toks@\@xp\@xp\@xp{\csname the#1\endcsname}%
  \@xp\xdef\csname the#1\endcsname
    {\@xp\@nx\csname the#2\endcsname
     .\the\toks@}}}}

\newcommand{\resetfootnoterule} {
  \renewcommand\footnoterule{%
  \kern-3\p@
  \hrule\@width.4\columnwidth
  \kern2.6\p@}
}

\makeatother
\numberwithin{equation}{section}

\theoremstyle{definition}
\newtheorem{remark}{Remark}

\begin{document}

\thispagestyle{empty}
\renewcommand{\ehkol}{E.V.\ Ferapontov and A.M.\ Grundland}
\renewcommand{\ohkol}{Different Analytic Descriptions 
of Mean Curvature Surfaces}

\begin{flushleft}
\footnotesize \sf
Journal of Nonlinear Mathematical Physics \qquad 2000, V.7, N~1,
\pageref{fer_grund_fp}--\pageref{fer_grund_lp}.
\hfill {\sc Letter}
\end{flushleft}

\vspace{-5mm}

\renewcommand{\footnoterule}{}
{\renewcommand{\thefootnote}{}
 \footnotetext{\prava{E.V.\ Ferapontov and A.M.\ Grundland}}}

\name{Links Between Different Analytic Descriptions \\
of Constant Mean Curvature Surfaces} \label{fer_grund_fp}

\Author{E.V.\ FERAPONTOV~$^{\dag}$ and A.M.\ GRUNDLAND~$^{\ddag}$}

\Adress{$^{\dag}$ Department of Mathematical Sciences,
Loughborough University, Loughborough, \\
~~Leicestershire LE11 3TU, United Kingdom \\
~~e-mail: E.V.Ferapontov@lboro.ac.uk \\[2mm]
~~Centre for Nonlinear Studies,
Landau Institute for Theoretical Physics, \\
~~Academy of Sciences of Russia, Kosygina 2,
117940 Moscow, GSP-1, Russia \\
~~e-mail: fer@landau.ac.ru \\[2mm]
$^{\ddag}$ Centre de Recherches Mathematiques,
Universit\'{e} de Montr\'{e}al,
C.P. 6128 Succ. \\
~~Centre-Ville,
2910 Chemin de la Tour, Montr\'{e}al, Qu\'{e}bec
H3C 3J7 Canada \\
~~e-mail: grundlan@crm.umontreal.ca}

\Date{Received September 6, 1999; Revised October 5, 1999; 
Accepted November 5, 1999}

\begin{abstract}
\noindent
Transformations between different analytic descriptions
of constant mean curvature (CMC) surfaces are established.
In particular, it is demonstrated that the system
\[
\begin{split}
&\partial \psi_{1} = (|\psi_{1}|^{2} + |\psi_{2}|^{2}) \psi_{2} \\
&\bar{\partial} \psi_{2} =- (|\psi_{1}|^{2} + |\psi_{2}|^{2}) \psi_{1}
\end{split}
\]
descriptive of CMC surfaces within the framework of the generalized
Weierstrass representation, decouples into a direct sum of
the elliptic Sh-Gordon and Laplace equations. Connections of this
system with the sigma model equations are established. It is
pointed out, that the instanton solutions correspond to
different Weierstrass parametrizations of the standard sphere
$S^{2} \subset E^{3}$.
\end{abstract}

\section{Introduction}

We investigate the system
\begin{equation} \label{eq1}
\begin{split}
&\partial \psi_{1} = (|\psi_{1}|^{2} + |\psi_{2}|^{2}) \psi_{2} \\
&\bar{\partial} \psi_{2} = - (|\psi_{1}|^{2} + |\psi_{2}|^{2} ) \psi_{1}
\end{split}
\end{equation}
which has been derived in [1] and governs constant mean
curvature (CMC) surfaces  in the conformal parametrization
$z, \bar{z}$ $(\partial = \partial_{z}, \bar{\partial} =
\partial_{\bar{z}}$). This system was subsequently
discussed in [2,4,6].
In this paper, we demonstrate that system (\ref{eq1})
can be decoupled into a direct sum of the
elliptic Sh-Gordon and Laplace equations.
Firstly, we change from $\psi_{1}$, $\psi_{2}$ to the
new dependent variables $Q,R$
\[
Q= 2 (\psi_{2} \partial \bar{\psi}_{1} - \bar{\psi}_{1}
\partial \psi_{2} ),
\quad
R = 2 (|\psi_{1}|^{2} + |\psi_{2}|^{2})^{2}/ |Q|.
\]
Secondly, we introduce new independent variables $\eta$, $\bar{\eta}$
according to
\[
d \, \eta = \sqrt{Q} \, dz,
\quad
d \bar{\eta} = \sqrt{\bar{Q}} \, d \bar{z}
\]
(these formulae are correct since $Q$ is holomorphic).
In the new variables $Q,R,\eta,\bar{\eta}$, system (\ref{eq1})
assumes the decoupled form
\begin{gather*}
(\ln R)_{\eta \bar{\eta}} = \frac{1}{R} - R, \\
Q_{\bar{\eta}} = \bar{Q}_{\eta} = 0,
\end{gather*}
which is a direct sum of the elliptic sh-Gordon and
Laplace equations.
This transformation is an immediate corollary of the known
properties of CMC surfaces.
Connection of system (\ref{eq1})
with the sigma-model equations
\[
\partial \bar{\partial} \rho - \frac{2 \bar{\rho}}{1 + |\rho|^{2}}
\partial \rho \bar{\partial} \rho = 0,
\]
($\rho =i \bar{\psi}_{1}/ {\psi}_{2}$) is also discussed. In terms
of $\rho$, our transformation adopts the form:
\[
Q = \frac{2 \partial \rho \, \partial \bar{\rho}}{(1+ |\rho|^{2})^{2}},
\quad
R = \left| \frac{ \partial \bar{\rho}}
{\partial \rho} \right|.
\]

\section{Generalized Weierstrass representation of surfaces in
$\mathbb R^{3}$}

Following the results of [1], with any solution $\psi_{1}$, $\psi_{2}$
of
the Dirac equations
\begin{equation} \label{eq2}
\partial \psi_{1} = p \psi_{2}, \quad
\bar{\partial} \psi_{2} = - p \psi_{1},
\end{equation}
($p(z, \bar{z})$ is a real potential),
we associate a surface $M^{2} \subset E^{3}$ with the radius-vector
${\bf r}(z , \bar{z})$ defined by the formulae
\begin{gather*}
\partial {\bf r} = {\big(} i (\psi_{2}^{2} + \bar{\psi}_{1}^{2}), \quad
\bar{\psi}_{1}^{2} - \psi_{2}^{2}, \quad -2 \psi_{2} \bar{\psi}_{1}
{\big)}, \\
\bar{\partial} {\bf r} = {\big(} -i (\bar{\psi}_{2}^{2} + \psi_{1}^{2}),
\quad
\psi_{1}^{2} - \bar{\psi}_{2}^{2}, \quad -2 \psi_{1} \bar{\psi}_{2}
{\big)}.
\end{gather*}
The latter are compatible by virtue of (\ref{eq2}). The unit normal ${\bf n}$ of
the surface $M^{2}$ can be calculated according to the
formula ${\bf n} = \frac{1}{i} \partial {\bf r} \times \bar{\partial}
{\bf r}/
 |\partial {\bf r} \times \bar{\partial} {\bf r}|$ and is represented
as follows:
\begin{gather}
{\bf n} = \frac{1}{q} \big( i (\bar{\psi}_{1} \bar{\psi}_{2} - \psi_{1}
\psi_{2}
), \quad
\bar{\psi}_{1} \bar{\psi}_{2} + \psi_{1} \psi_{2}, \quad
\psi_{1} \bar{\psi}_{1} - \psi_{2} \bar{\psi}_{2} {\big)}, \label{eq3} \\
q= |\psi_{1}|^{2} + |\psi_{2} |^{2}. \notag
\end{gather}
One can verify directly that the scalar products in $E^{3}$ are
$({\bf n},{\bf n})=1$, $({\bf n}, \partial {\bf r})=
({\bf n}, \bar{\partial} {\bf r}) = 0$.
Equations of motion of the complex frame $\partial {\bf r}$,
$\bar{\partial} {\bf r}$,
${\bf n}$ are of the form
\begin{equation} \label{eq4}
\begin{split}
&\partial \left(  \begin{array}{c}
\partial {\bf r}  \\
\bar{\partial} {\bf r}  \\
{\bf n}    \\
\end{array}
    \right)  =
         \left(  \begin{array}{ccc}
2 \displaystyle \frac{\partial q}{q} &  0  &   Q   \\
0       &     0    &  2 H q^{2}      \\
-H      &  - \displaystyle \frac{Q}{2 q^{2}}  &  0  \\
\end{array}   \right)
         \left(  \begin{array}{c}
\partial {\bf r}    \\
\bar{\partial} {\bf r}  \\
{\bf n}             \\
\end{array}   \right)
\\[1ex]
&\bar{\partial} \left( \begin{array}{c}
\partial {\bf r}  \\
\bar{\partial} {\bf r}   \\
{\bf n}     \\
\end{array}
   \right)  =
         \left(  \begin{array}{ccc}
0  &    0 &   2 H q^{2}    \\
0  &  2 \displaystyle \frac{\bar{\partial} q}{q}  &  \bar{Q}   \\
- \displaystyle \frac{\bar{Q}}{2 q^{2}} & -H  &   0  \\
\end{array}   \right)
         \left(  \begin{array}{c}
\partial {\bf r}    \\
\bar{\partial} {\bf r}  \\
{\bf n}    \\
\end{array}   \right)
\end{split}
\end{equation}
where $Q= 2 (\psi_{2} \partial \bar{\psi}_{1} - \bar{\psi}_{1}
\partial \psi_{2})$ and $H = p/q$ is the mean curvature.
Formulae (\ref{eq4}) are compatible with the scalar products
\[
({\bf n},{\bf n})=1,  \quad (\partial {\bf r}, \bar{\partial} {\bf r}) =
2 q^{2}
\]
(all other scalar products being equal to zero). Using (\ref{eq4}),
one can derive the following useful equation for the unit
normal ${\bf n}$:
\begin{equation} \label{eq5}
\partial \bar{\partial} {\bf n} +
(\partial {\bf n} , \bar{\partial} {\bf n}) {\bf n} +
\bar{\partial} H \, \partial {\bf r} +
\partial H \, \bar{\partial} {\bf r} = 0.
\end{equation}

The first fundamental form $I= (d {\bf r}, d {\bf r})$ and the second
fundamental form $II= (d^{2} {\bf r}, {\bf n})$ of the surface
$M^{2}$ are given by
\begin{equation} \label{eq6}
\begin{split}
&I = 4 q^{2} \, dz \, d \bar{z}, \\
&II = Q \, dz^{2} +  4 H q^{2} dz \, d \bar{z} + \bar{Q} \, d
\bar{z}^{2}.
\end{split}
\end{equation}
The quantity $Q \, dz^{2}$ is called the Hopf differential.
The real potential $p$ and the spinors $\psi_1, \psi_2$ satisfying (\ref{eq1}) can be
viewed as the "generalized Weierstrass data" of the surface $M^2$.
The corresponding Gauss-Codazzi equations
which are the compatibility conditions for (\ref{eq4}), are of the form
\begin{equation} \label{eq7}
\begin{split}
&\partial \bar{\partial} ( \ln q^{2}) = \frac{1}{2} \frac{Q
\bar{Q}}{q^{2}}
-2 H^{2} q^{2}, \\
&\bar{\partial} Q = 2  q^{2} \partial H , \\
&\partial \bar{Q} = 2  q^{2} \bar{\partial} H.
\end{split}
\end{equation}
In fact, equations (\ref{eq7}) are a direct differential consequence of (\ref{eq2}), as
can be
checked by a straightforward calculation.
Gauss-Codazzi equations of surfaces in
conformal para\-metr\-iza\-tion $z, \bar{z}$ have been discussed in
[5]. We recall also that the Gaussian
curvature $K$ of the surface $M^{2}$ can be calculated as follows:
\begin{equation} \label{eq8}
K = - \frac{1}{q^{2}} \partial \bar{\partial} ( \ln q).
\end{equation}

The unit normal ${\bf n} = (n_{1}, n_{2}, n_{3})$, given by (\ref{eq3}), maps a
surface $M^{2}$ onto the unit sphere $S^{2}$. Combining this map with
the
stereographic projection, we obtain a map $\rho$ of the surface
$M^{2}$ onto the complex plane, called the complex Gauss map.
In our notation, $\rho$ assumes the form
\begin{equation} \label{eq9}
\rho = \frac{n_{1} + i n_{2}}{1 - n_{3}} = i
\frac{\bar{\psi}_{1}}{\psi_{2}}.
\end{equation}
According to the results of [3], Gauss map $\rho$
satisfies the nonlinear equation
\begin{equation} \label{eq10}
 (\partial \bar{\partial} \rho - \frac{2 \bar{\rho}}{1 + |\rho|^{2}}
\partial \rho \, \bar{\partial} \rho) H= \partial H \, \bar{\partial}
\rho,
\end{equation}
which formally can be viewed as a differential consequence of (\ref{eq2}).
In terms of $\rho$ and $H$, the initial data $\psi_{1}$, $\psi_{2}$ and
$p$ assume the form
\[
\psi_{1} = \frac{ \bar{\rho}}{\sqrt{H}} \frac{\sqrt{i \bar{\partial}
\rho}}
{1 + |\rho|^{2}},
\qquad
\psi_{2} = \frac{1}{\sqrt{H}} \frac{\sqrt{i \partial \bar{\rho}}}
{1 + |\rho|^{2}},
\qquad
p = \frac{|\partial \bar{\rho}|}{1+ |\rho|^{2}}
\]
while the expressions for $q$ and $Q$ take the form
\[
q= \frac{1}{H} \frac{\bar{\partial} \rho \partial \bar{\rho}}
{1 + |\rho|^{2}},
\quad
Q = \frac{2}{H} \frac{\partial \rho \, \partial \bar{\rho}}{(1+
|\rho|^{2})^{2}}.
\]
\begin{remark}
Linear problem (\ref{eq4}) can be rewritten in terms of
$\psi_{1}$, $\psi_{2}$ as follows. First of all, we
point out that
\[
\partial q = \psi_{1} \partial \bar{\psi}_{1} + \bar{\psi}_{2}
\partial \psi_{2}.
\]
Combining this equation with the definition of $Q$:
\[
\frac{1}{2} Q = \psi_{2} \partial  \bar{\psi}_{1} - \bar{\psi}_{1}
\partial \psi_{2},
\]
and solving these two equations for $\partial \bar{\psi}_{1}$,
$\partial \psi_{2}$, we can ``close'' system (\ref{eq2}) as follows:
\begin{equation} \label{eq11}
\begin{split}
&\partial \left(
\begin{array}{c}
\psi_{1}   \\
\psi_{2}   \\
\end{array}   \right) =
\left(
\begin{array}{cc}
0   &  q H    \\
- \displaystyle \frac{Q}{2q} & \displaystyle \frac{\partial q}{q}  \\
\end{array}   \right)
\left(
\begin{array}{c}
\psi_{1}   \\
\psi_{2}   \\
\end{array}  \right),
\\[1ex]
&\bar{\partial} \left(
\begin{array}{c}
\psi_{1}   \\
\psi_{2}   \\
\end{array}   \right) =
\left(
\begin{array}{cc}
\displaystyle \frac{\bar{\partial} q}{q} & \displaystyle
\frac{\bar{Q}}{2q}  \\
-qH   &   0   \\
\end{array}    \right)
\left(
\begin{array}{c}
\psi_{1}   \\
\psi_{2}   \\
\end{array}   \right).
\end{split}
\end{equation}
The compatibility conditions for system (\ref{eq11}) coincide with (\ref{eq7}).
We point out that the $2 \times 2$ matrix approach to surfaces in
$E^{3}$
has been extensively developed in [5]. From the point of view of the
theory of integrable systems, linear system (\ref{eq4}) can be regarded
as the squared eigenfunction equations corresponding to (\ref{eq11})
(indeed, $\partial {\bf r}$, $\bar{\partial} {\bf r}$ and ${\bf n}$
are quadratic
expressions in $\psi_{1}$ and $\psi_{2}$).
\end{remark}

\section{CMC-1 surfaces}

The class of CMC-1 surfaces is characterized by the constraint $H=1$ or,
equivalently, $p=q$. Introducing this ansatz in (\ref{eq2}), we arrive at
the nonlinear system (\ref{eq1})
\begin{gather*}
\partial \psi_{1} = (|\psi_{1}|^{2} + |\psi_{2}|^{2}) \psi_{2}, \\
\bar{\partial} \psi_{2} = - (|\psi_{1}|^{2} + |\psi_{2}|^{2}) \psi_{1},
\end{gather*}
which is the main subject of our study. According to the
previous section, system (\ref{eq1}) is equivalent to
\begin{equation} \label{eq12}
\begin{split}
&\partial \bar{\partial} \ln q^{2} =
\displaystyle \frac{1}{2} \frac{Q \bar{Q}}{q^{2}}
- 2 q^{2},  \\
&\bar{\partial} Q = \partial \bar{Q} = 0,
\end{split}
\end{equation}
where $q= |\psi_{1}|^{2} + |\psi_{2}|^{2}$, $Q= 2(\psi_{2}
\partial \bar{\psi}_{1} - \bar{\psi}_{1} \partial \psi_{2})$.
Thus, for CMC surfaces the Hopf differential $Q \, dz^{2}$ is
holomorphic.
Applying to system (\ref{eq12}) the reciprocal transformation
\[
d \eta = \sqrt{Q} \, dz , \quad
d \bar{\eta} = \sqrt{\bar{Q}} \, d \bar{z}
\]
(that is, changing from $z,\bar{z}$ to the new independent
variables $\eta, \bar{\eta}$ which are correctly defined
in view of the holomorphicity of $Q$) and introducing
\[
R= \frac{2 q^{2}}{|Q|},
\]
we transform system (\ref{eq12}) into the decoupled form
\begin{equation} \label{eq13}
\begin{split}
&(\ln R)_{\eta \bar{\eta}} = \displaystyle \frac{1}{R} - R, \\
&Q_{\bar{\eta}} = \bar{Q}_{\eta} = 0. 
\end{split}
\end{equation}
This result provides the rationale for the change of variables which
we introduce in Section~1.

\begin{remark}
System (\ref{eq1}) is invariant under the $SU(2)$-symmetry
\begin{equation} \label{eq14}
\zeta_{1} = \alpha \psi_{1} + \beta \bar{\psi}_{2}, \quad
\zeta_{2} = - \beta \bar{\psi}_{1} + \alpha \psi_{2},
\end{equation}
where $\alpha, \beta$ are complex constants subject to
the constraint $\alpha \bar{\alpha} + \beta \bar{\beta} = 1$.
One can check directly, that the quantities $Q$ and $R$ are
invariant under transformations (\ref{eq14}), so that the surfaces,
corresponding to $(\psi_{1}, \psi_{2})$ and $(\zeta_{1}, \zeta_{2})$,
have coincident fundamental forms. Thus, they are  identical
up to a rigid motion in $E^{3}$. The passage from
$\psi_{1}, \psi_{2}$ to $Q,R$ can thus be viewed as a passage to
the differential invariants of the point symmetry group (\ref{eq14}).
\end{remark}

\begin{remark}
For CMC-1 surfaces, system (\ref{eq11}) allows the
introduction of a spectral parameter
\begin{equation} \label{eq15}
\begin{split}
&\partial \left(
\begin{array}{c}
\psi_{1}  \\
\psi_{2}   \\
\end{array} \right)
= \left(
\begin{array}{cc}
0   & q   \\
\displaystyle - \lambda \frac{Q}{2q}  &  \displaystyle \frac{\partial
q}{q} \\
\end{array} \right)
\left(
\begin{array}{c}
\psi_{1}  \\
\psi_{2}  \\
\end{array}  \right)
\\[1ex]
&\bar{\partial} \left(
\begin{array}{c}
\psi_{1}  \\
\psi_{2}  \\
\end{array}   \right)
= \left(
\begin{array}{cc}
\displaystyle \frac{\bar{\partial} q}{q} &
\displaystyle \frac{1}{\lambda} \frac{\bar{Q}}{2q}  \\
-q   &   0     \\
\end{array}  \right)
\left(
\begin{array}{c}
\psi_{1}  \\
\psi_{2}  \\
\end{array}   \right)
\end{split}
\end{equation}
where $\lambda$ is a unitary constant, $|\lambda|=1$.
The gauge transformation
\[
\tilde{\psi}_{1} = \psi_{1},  \quad
\tilde{\psi}_{2} = q^{-1} \, \psi_{2}
\]
reduces linear spectral problem (\ref{eq15}) to the $SL(2)$ form
\begin{equation} \label{eq16}
\begin{split}
&\partial       \left(
\begin{array}{c}
\tilde{\psi}_{1}  \\
\tilde{\psi}_{2}  \\
\end{array}    \right)  =
\left(
\begin{array}{cc}
0     &   q^{2}    \\
      &            \\
- \displaystyle  \lambda \frac{Q}{2 q^{2}}  &  0   \\
\end{array}    \right)
\left(
\begin{array}{c}
\tilde{\psi}_{1}   \\
\tilde{\psi}_{2}   \\
\end{array}    \right)
\\[1ex]
&\bar{\partial}  \left(
\begin{array}{c}
\tilde{\psi}_{1}   \\
\tilde{\psi}_{2}   \\
\end{array}    \right)  =
\left(
\begin{array}{cc}
\displaystyle \frac{\bar{\partial} q}{q}  &
\displaystyle \frac{\bar{Q}}{2 \lambda}  \\
    &             \\
-1  &  - \displaystyle \frac{\bar{\partial} q}{q}    \\
\end{array}    \right)
\left(
\begin{array}{c}
\tilde{\psi}_{1}   \\
\tilde{\psi}_{2}   \\
\end{array}    \right).
\end{split}
\end{equation}
The compatibility conditions for both systems (\ref{eq15}) and (\ref{eq16}) coincide with (\ref{eq12}).
From the linear system (\ref{eq16}) the radius-vector $\bf r$ can be recovered
via the so-called Sym formula: we refer to [5] and [9] for the further
discussion of the Sym approach.

Linear system (\ref{eq16}) can be readily rewritten in terms of $\psi_{1}$,
$\psi_{2}$. Indeed, observing that (\ref{eq16}) implies
\[
\partial \ln q = \partial \ln (|\psi_{1}|^{2}  + |\psi_{2} |^{2}),
 \quad
\bar{\partial} \ln q = \bar{\partial} \ln (|\psi_{1}|^{2} +
|\psi_{2}|^{2} ),
\]
we can take $q= c( |\psi_{1}|^{2}  + |\psi_{2}|^{2})$,
$c \in \mathbb C$, which, upon substitution in (\ref{eq16}), produces system
(\ref{eq1}). Thus
transformation from (\ref{eq12}) to (\ref{eq1}) consists of rewriting (\ref{eq16}) in terms of
the
$\psi$. Representations in terms of $\psi$ are called
eigenfunction equations, and are fundamental in soliton theory
-- see eg [7]. The Lax pair for system (\ref{eq1}) is of the form [6]
\begin{gather*}
\partial \left(
\begin{array}{c}
\phi_{1}  \\
\phi_{2}  \\
\end{array}   \right)
= \frac{2}{\mu +1} \left(
\begin{array}{cc}
- \displaystyle \bar{\psi}_{1} \psi_{2} + \frac{Q}{2 q^{2}} \psi_{1}
\bar{\psi}_{2} &
- \bar{\psi}_{1}^{2} - \displaystyle \frac{Q}{2 q^{2}}
\bar{\psi}_{2}^{2}   \\
             &     \\
\displaystyle \psi_{2}^{2} + \frac{Q}{2 q^{2}} \psi_{1}^{2}    &
\bar{\psi}_{1} \psi_{2} - \displaystyle \frac{Q}{2 q^{2}} \psi_{1}
\bar{\psi}_{2} \\
\end{array}    \right)
\left(
\begin{array}{c}
\phi_{1}  \\
\phi_{2}  \\
\end{array}   \right)
\\[1ex]
\bar{\partial} \left(
\begin{array}{c}
\phi_{1}   \\
\phi_{2}   \\
\end{array}   \right)
= \frac{2}{\mu - 1} \left(
\begin{array}{cc}
- \displaystyle \psi_{1} \bar{\psi}_{2} + \frac{\bar{Q}}{2 q^{2}}
\bar{\psi}_{1} \psi_{2} &
\bar{\psi}_{2}^{2} + \displaystyle \frac{\bar{Q}}{2 q^{2}}
\bar{\psi}_{1}^{2}    \\
      &                        \\
- \displaystyle \psi_{1}^{2} - \frac{\bar{Q}}{2 q^{2}} \psi_{2}^{2}
&
\psi_{1} \bar{\psi}_{2} - \displaystyle \frac{\bar{Q}}{2 q^{2}}
\bar{\psi}_{1} \psi_{2} \\
\end{array}    \right)
\left(
\begin{array}{c}
\phi_{1}    \\
\phi_{2}    \\
\end{array}    \right).
\end{gather*}
It is interesting to note that the compatibility conditions for this
linear system, which is of the first order in the derivatives of
$\psi$, give us exactly system (\ref{eq1}). The latter is also of the
first order in $\psi$.
\end{remark}

\section{CMC surfaces and sigma model equations}

For CMC surfaces, equations (\ref{eq10}) imply the nonlinear sigma model
\begin{equation} \label{eq17}
\partial \bar{\partial} \rho - \frac{2 \bar{\rho}}{1 +|\rho|^{2}}
\partial \rho \, \bar{\partial} \rho = 0,
\end{equation}
descriptions of the stationary two-dimensional $SU(2)$ magnet.
The transformation of the system (\ref{eq17}) into the decoupled form (\ref{eq13}) now
assumes the form
\begin{equation} \label{eq18}
\begin{split}
&Q = 2 \displaystyle \frac{\partial \rho \, \partial
\bar{\rho}}{(1+|\rho|^{2})^{2}},
\qquad
R= |\displaystyle \frac{\partial \bar{\rho}}
{\partial \rho }|,  \\
&d \eta = \sqrt{Q} \, dz ,  \qquad
d \bar{\eta} = \sqrt{\bar{Q}} \, d \bar{z}.  \\
\end{split}
\end{equation}
In terms of the unit normal vector ${\bf n}$, equation (\ref{eq5})
adopts the form of the $SO(3)$ sigma model
\begin{equation} \label{eq19}
\bar{\partial} \partial {\bf n} + ( \partial {\bf n}, \bar{\partial}
{\bf n} )
{\bf n} = 0,
\quad
({\bf n}, {\bf n}) = 1.
\end{equation}
Formula (\ref{eq9}) establishes a link between sigma models (\ref{eq17}) and (\ref{eq19}).
The topological charge
\[
\frac{1}{4 \pi} \int \int ({\bf n}, [ \partial {\bf n} \times
\bar{\partial} {\bf n}]) \,
dz \wedge d \bar{z}
\]
can be written as
\[
\frac{1}{2 \pi i} \int \int \partial \bar{\partial} \ln \, q \:
dz \wedge d \bar{z}
\]
which, in view of (\ref{eq8}), is the topologically invariant
integral curvature of the surface
$M^{2}$.

Instanton solutions of system (\ref{eq19}) are specified by the ansatz
\[
\partial {\bf n} = \pm i {\bf n} \times \partial {\bf n},
\quad
\bar{\partial} {\bf n} = \mp i {\bf n} \times \bar{\partial} {\bf n},
\]
which, after a simple calculation, implies $Q=0$. Solutions of system
(\ref{eq1})
specified by a constraint $Q=0$, can be represented in the form
\begin{equation} \label{eq20}
\psi_{1} = \frac{\rho \sqrt{\bar{\partial} \bar{\rho}}}{1 + |\rho|^{2}},
\quad
\psi_{2} = \frac{\sqrt{\partial \rho}}{1 + |\rho|^{2}},
\quad
p = \frac{|\partial \rho|}{1+ |\rho|^{2}},
\end{equation}
where $\rho(z)$ is an arbitrary holomorphic function.
In the case when the energy
\[
E = \int \int \frac{\partial \rho \, \bar{\partial} \rho}
{1 + |\rho|^{2}} \, dz \wedge d \bar{z}
\]
is finite, the function $\rho(z)$ is rational in $z$ [8].
\par
Geometrically, instanton solutions (\ref{eq20}) parametrize the
standardly embedded sphere $S^{2} \subset E^{3}$.
This can be readily seen from formulae (\ref{eq6}), which, in case $Q=0$, imply
the
proportionality of fundamental forms I and II. This example shows that
different
Weierstrass data $(\psi_{1}, \psi_{2}, p)$ can correspond to
different parametrizations of one and the same surface $M^{2} \subset
E^{3}$.

Introducing the two-component complex vector
\[
{\bf N} = ( \frac{\psi_{1}}{\sqrt{q}}, \frac{\bar{\psi}_{2}}{\sqrt{q}}),
\quad
q = |\psi_{1}|^{2} + |\psi_{2}|^{2},
\]
one can check that ${\bf N}$ satisfies the equations of the $\mathbb C
P^{1}$
sigma model
\begin{equation} \label{eq21}
\begin{split}
&( {\bf N}, \bar{\bf N} ) = 1,  \\
&\partial \bar{\partial} {\bf N} = ( \bar{\bf N}, \bar{\partial} {\bf N})
\, \partial {\bf N} + ( \bar{\bf N}, \partial {\bf N}) \, \bar{\partial}
{\bf N}
- k {\bf N},
\end{split}
\end{equation}
where
\[
k = -2 ( \bar{ \bf N}, \partial {\bf N}) ( {\bf N}, \bar{\partial}
\bar{\bf N}
) + \frac{1}{2} ( \partial {\bf N}, \bar{\partial} \bar{\bf N}) +
\frac{1}{2}( \bar{\partial} {\bf N}, \partial \bar{\bf N}).
\]
Equations (\ref{eq21}) are associated with the Lagrangian
\[
L = \int \int \{ (\bar{\partial} {\bf N}, \partial \bar{\bf N}) +
( \partial {\bf N}, \bar{\partial} \bar{\bf N}) +
2 ( \bar{\bf N}, \partial {\bf N}) (\bar{\bf N},
\bar{\partial} {\bf N}) - 2 k
[ ( {\bf N}, \bar{\bf N}) -1] \} \, dz \, d \bar{z},
\]
where $k$ is the Lagrange multiplier.

\subsection*{Acknowledgments} 
One of the authors (E.V.F.) would like to
thank Centre de Recherches Math\'{e}matiques,
Universit\'{e} de Montr\'{e}al for generous hospitality
and financial support during the period when this investigation
was performed. This work was supported by research grant from NSERC
of Canada and the Fonds FCAR du Gouvernment du Qu\'{e}bec.
We would like to thank Professor B. G. Konopelchenko for drawing our
attention to this subject and for helpful discussions. We would
like to thank the referee for useful comments.

\label{fer_grund_lp}


\begin{thebibliography}{9}

\footnotesize

\bibitem{1}
Konopelchenko B.G., Induced Surfaces and Their Integrable
Dynamics, {\it Stud.\ Appl.\ Math.}, 1996, V.96, 9--51.
\bibitem{2}
Konopelchenko B.G.\ and Taimanov I.A.,
Constant Mean Curvature Surfaces via an Integrable
Dynamical System, {\it J.\ Phys A: Math.\ Gen.}, 1998, V.29, 1261--1265. 
\bibitem{3}
Kenmotsu K., Weierstrass Formula for Surfaces of
Prescribed Mean Curvature, {\it Math.\ Ann.}, 1979, V.245, 89--99.
\bibitem{4}
Carrol R.\ and Konopelchenko B.G., Generalized
Weierstrass-Enneper Inducing, Conformal Immersions and Gravity,
{\it Int.\ J.\ Modern Phys.\ A}, 1996, V.11, N~7, 1183--1216.
\bibitem{5}
Bobenko A.I., Surfaces in Terms of $2 \times 2$ Matrices.
Old and New Integrable Cases, in Harmonic Maps and Integrable
Systems, Editors A.\ Fordy and J.\ Wood, Vieweg, 1994, 83--127.
\bibitem{6}
Bracken P., Grundland A.M.\ and Martina L.,
The Weierstrass-Enneper System for Constant Mean Curvature Surfaces and
the Completely Integrable Sigma Model, {\it J.\ Math.\ Phys.}, 1999, V.40, N~7,
3379--3403.
\bibitem{7}
Konopelchenko B.G., Soliton Eigenfunction Equations: The IST
Integrability and Some Properties, {\it Rev.\ in Math.\ Phys.}, 1990, V.2,
399--466.
\bibitem{8}
Zakrzewski W., Low Dimensional Sigma-Models, Hilger, London,
1989.
\bibitem{9}
Sym A., Soliton Surfaces and Their Applications, in Geometric Aspects
of the Einstein Equations and Integrable Systems, Lecture Notes in Phys, 
239, Springer-Verlag, Berlin, 1985, 154--231.

\end{thebibliography}
\end{document}